\documentstyle[amsfonts,12pt]{article}
\begin{document}

\title{Stochastic solution of a nonlinear fractional differential equation}
\author{F. Cipriano\thanks{%
GFM and FCT-Universidade Nova de Lisboa, Complexo Interdisciplinar, Av. Gama
Pinto, 2 - 1649-003 Lisboa (Portugal), e-mail: cipriano@cii.fc.ul.pt}, H.
Ouerdiane\thanks{%
D\'{e}partment de Math\'{e}matiques, Facult\'{e} des Sciences,
Universit\'{e} de Tunis El Manar, Campus Universitaire, 1060 Tunis. e-mail:
Habib.Ouerdiane@fst.rnu.tn} and R. Vilela Mendes\thanks{%
CMAF, Complexo Interdisciplinar, Universidade de Lisboa, Av. Gama Pinto, 2 -
1649-003 Lisboa (Portugal), e-mail: vilela@cii.fc.ul.pt,
http://label2.ist.utl.pt/vilela/} \thanks{%
Centro de Fus\~{a}o Nuclear - EURATOM/IST Association, Instituto Superior
T\'{e}cnico, Av. Rovisco Pais 1, 1049-001 Lisboa, Portugal}}
\date{}
\maketitle

\begin{abstract}
A stochastic solution is constructed for a fractional generalization of the
KPP (Kolmogorov, Petrovskii, Piskunov) equation. The solution uses a
fractional generalization of the branching exponential process and
propagation processes which are spectral integrals of Levy processes.
\end{abstract}

\section{Introduction: The notion of stochastic solution}

The solutions of linear elliptic and parabolic equations, both with Cauchy
and Dirichlet boundary conditions, have a probabilistic interpretation.
These are classical results which may be traced back to the work of Courant,
Friedrichs and Lewy \cite{Courant} in the 1920's and became a standard tool
in potential theory\cite{Getoor} \cite{Bass2}. For example, for the heat
equation

\begin{equation}
\partial _{t}u(t,x)=\frac{1}{2}\frac{\partial ^{2}}{\partial x^{2}}%
u(t,x)\qquad \textnormal{with}\qquad u(0,x)=f(x)  \label{1.1}
\end{equation}
the solution may be written either as 
\begin{equation}
u\left( t,x\right) =\frac{1}{2\sqrt{\pi }}\int \frac{1}{\sqrt{t}}\exp \left(
-\frac{\left( x-y\right) ^{2}}{4t}\right) f\left( y\right) dy  \label{1.2}
\end{equation}
or as 
\begin{equation}
u(t,x)={\Bbb E}_{x}f(X_{t})  \label{1.3}
\end{equation}
${\Bbb E}_{x}$ meaning the expectation value, starting from $x$, of the
process 
\[
dX_{t}=dW_{t} 
\]
$W_{t}$ being the Wiener process.

Eq.(\ref{1.1}) is a {\it specification} of a problem whereas (\ref{1.2}) and
(\ref{1.3}) are {\it solutions} in the sense that they both provide
algorithmic means to construct a function satisfying the specification. An
important condition for (\ref{1.2}) and (\ref{1.3}) to be considered as
solutions is the fact that the algorithmic tools are independent of the
particular solution, in the first case an integration procedure and in the
second the simulation of a solution-independent process. This should be
contrasted with stochastic processes constructed from a given particular
solution, as has been done for example for the Boltzman equation\cite{Graham}%
.

In contrast with the linear problems, for nonlinear partial differential
equations, explicit solutions in terms of elementary functions or integrals
are only known in very particular cases. However, if a solution-independent
stochastic process is found that (for arbitrary initial conditions)
generates the solution in the sense of Eq.(\ref{1.3}), a stochastic solution
is obtained. In this way the set of equations for which exact solutions are
known would be considerably extended.

The stochastic representations recently constructed for the Navier-Stokes%
\cite{Jan} \cite{Waymire} \cite{Bhatta1} \cite{Ossiander} and the
Vlasov-Poisson equations\cite{Vilela1} \cite{Vilela2} define
solution-independent processes for which the mean values of some functionals
are solutions to these equations. Therefore, they are exact {\bf stochastic
solutions}.

In the stochastic solutions one deals with a process that starts from the
point where the solution is to be found, a functional being then computed
along the whole sample path or until it reaches a boundary. In all cases one
needs to average over many independent sample paths to obtain a expectation
value of the functional. The localized and parallelizable nature of the
solution construction is clear. Provided some differentiability conditions
are satisfied, the process also handles equally well simple or very complex
boundary conditions.

Stochastic solutions also provide an intuitive characterization of the
physical phenomena, relating nonlinear interactions with cascading
processes. By the study of exit times from a domain they also sometimes
provide access to quantities that cannot be obtained by perturbative methods%
\cite{VilelaZeit}.

One way to construct stochastic solutions is based on a probabilistic
interpretation of the Picard series. The differential equation is written as
an integral equation which is rearranged in a such a way that the
coefficients of the successive terms in the Picard iteration obey a
normalization condition. The Picard iteration is then interpreted as an
evolution and branching process, the stochastic solution being equivalent to
importance sampling of the normalized Picard series. This method is used in
this paper to obtain a stochastic solution of a nonlinear partial
differential equation, which is a fractional version of the
Kolmogorov-Petrovskii-Piskunov (KPP)\ equation\cite{KPP}.

\section{A fractional nonlinear partial differential equation}

We consider the following equation

\begin{equation}
_{t}D_{*}^{\alpha }u\left( t,x\right) =\frac{1}{2}\,_{x}D_{\theta }^{\beta
}u\left( t,x\right) +u^{2}\left( t,x\right) -u\left( t,x\right)  \label{2.1}
\end{equation}
We use the same notations as in the study of the linear problem in \cite
{Mainardi1}. $_{t}D_{*}^{\alpha }$ is a Caputo derivative of order $\alpha $%
\begin{equation}
_{t}D_{*}^{\alpha }f\left( t\right) =\left\{ 
\begin{array}{lll}
\frac{1}{\Gamma \left( m-\alpha \right) }\int_{0}^{t}\frac{f^{(m)}\left(
\tau \right) d\tau }{\left( t-\tau \right) ^{\alpha +1-m}} &  & m-1<\alpha <m
\\ 
\frac{d^{m}}{dt^{m}}f\left( t\right) &  & \alpha =m
\end{array}
\right.  \label{2.2}
\end{equation}
$m$ integer. $_{x}D_{\theta }^{\beta }$ is a Riesz-Feller derivative defined
through its Fourier symbol by 
\begin{equation}
{\cal F}\left\{ _{x}D_{\theta }^{\beta }f\left( x\right) \right\} \left(
k\right) =-\psi _{\beta }^{\theta }\left( k\right) {\cal F}\left\{ f\left(
x\right) \right\} \left( k\right)  \label{2.3}
\end{equation}
with $\psi _{\beta }^{\theta }\left( k\right) =\left| k\right| ^{\beta
}e^{i\left( \textnormal{sign}k\right) \theta \pi /2}$.

Eq.(\ref{2.1}) is a fractional version of the KPP equation, studied by
probabilistic means by McKean\cite{McKean}. Physically it describes a
nonlinear diffusion with growing mass and in our fractional generalization
it would represent the same phenomenon taking into account memory effects in
time and long range correlations in space.

As outlined in the introduction, the first step towards a probabilistic
formulation is the rewriting of Eq.(\ref{2.1}) as an integral equation
including the initial conditions. For this purpose we take the Fourier
transform $\left( {\cal F}\right) $ in space and the Laplace transform $%
\left( {\cal L}\right) $ in time obtaining 
\begin{equation}
s^{\alpha }\widetilde{\widehat{u}}\left( s,k\right) \left( s,k\right)
=s^{\alpha -1}\widehat{u}\left( 0^{+},k\right) -\frac{1}{2}\psi _{\beta
}^{\theta }\left( k\right) \widetilde{\widehat{u}}\left( s,k\right) -%
\widetilde{\widehat{u}}\left( s,k\right) +\int_{0}^{\infty }dte^{-st}{\cal F}%
\left( u^{2}\left( t,x\right) \right)   \label{2.4}
\end{equation}
\newline
where 
\[
\widehat{u}\left( t,k\right) ={\cal F}\left( u\left( t,x\right) \right)
=\int_{-\infty }^{\infty }e^{ikx}u\left( t,x\right) dx
\]
\[
\widetilde{u}\left( s,x\right) ={\cal L}\left( u\left( t,x\right) \right)
=\int_{0}^{\infty }e^{-st}u\left( t,x\right) dt
\]
This equation holds for $0<\alpha \leq 1$ or for $0<\alpha \leq 2$ with $%
\frac{\partial }{dt}u\left( 0^{+},x\right) =0$. Solving for $\widetilde{%
\widehat{u}}\left( s,k\right) $ one obtains an integral equation 
\begin{equation}
\widetilde{\widehat{u}}\left( s,k\right) \,=\frac{s^{\alpha -1}}{s^{\alpha
}+1+\frac{1}{2}\psi _{\beta }^{\theta }\left( k\right) }\widehat{u}\left(
0^{+},k\right) +\int_{0}^{\infty }dt\frac{e^{-st}}{s^{\alpha }+1+\frac{1}{2}%
\psi _{\beta }^{\theta }\left( k\right) }{\cal F}\left( u^{2}\left(
t,x\right) \right)   \label{2.5}
\end{equation}
Taking the inverse Fourier and Laplace\cite{Trujillo} transforms 
\begin{eqnarray}
u\left( t,x\right)  &=&E_{\alpha ,1}\left( -t^{\alpha }\right) \int_{-\infty
}^{\infty }dy{\cal F}^{-1}\left( \frac{E_{\alpha ,1}\left( -\left( 1+\frac{1%
}{2}\psi _{\beta }^{\theta }\left( k\right) \right) t^{\alpha }\right) }{%
E_{\alpha ,1}\left( -t^{\alpha }\right) }\right) \left( x-y\right) u\left(
0^{+},y\right)   \nonumber \\
&&+\int_{0}^{t}d\tau \left( t-\tau \right) ^{\alpha -1}E_{\alpha ,\alpha
}\left( -\left( t-\tau \right) ^{\alpha }\right)   \nonumber \\
&&\int_{-\infty }^{\infty }dy{\cal F}^{-1}\left( \frac{E_{\alpha ,\alpha
}\left( -\left( 1+\frac{1}{2}\psi _{\beta }^{\theta }\left( k\right) \right)
\left( t-\tau \right) ^{\alpha }\right) }{E_{\alpha ,\alpha }\left( -\left(
t-\tau \right) ^{\alpha }\right) }\right) \left( x-y\right) u^{2}\left( \tau
,y\right)   \label{2.6}
\end{eqnarray}
$E_{\beta ,\rho }$ is the generalized Mittag-Leffler function 
\[
E_{\alpha ,\rho }\left( z\right) =\sum_{j=0}^{\infty }\frac{z^{j}}{\Gamma
\left( \alpha j+\rho \right) }
\]
We define the following propagation kernel 
\begin{equation}
G_{\alpha ,\rho }^{\beta }\left( t,x\right) ={\cal F}^{-1}\left( \frac{%
E_{\alpha .\rho }\left( -\left( 1+\frac{1}{2}\psi _{\beta }^{\theta }\left(
k\right) \right) t^{\alpha }\right) }{E_{\alpha ,\rho }\left( -t^{\alpha
}\right) }\right) \left( x\right)   \label{2.7}
\end{equation}
and, from the normalization relation, 
\[
E_{\alpha ,1}\left( -t^{\alpha }\right) +\int_{0}^{t}d\tau \left( t-\tau
\right) ^{\alpha -1}E_{\alpha ,\alpha }\left( -\left( t-\tau \right)
^{\alpha }\right) =1
\]
we may interpret $E_{\alpha ,1}\left( -t^{\alpha }\right) $ and $\left(
t-\tau \right) ^{\alpha -1}E_{\alpha ,\alpha }\left( -\left( t-\tau \right)
^{\alpha }\right) $, respectively as a survival probability up to time $t$
and as the probability density for the branching at time $\tau $ in a
branching process $B_{\alpha }$. It is a fractional generalization of an
exponential process. This provides a probabilistic sampling of the Picard
series obtained by iteration of Eq.(\ref{2.6}). The solution is therefore
obtained by the expectation of the exit values of the following process:

Starting at time zero, a particle lives according to the process $B_{\alpha }
$. At the branching time $\tau $ the initial particle dies and two new
particles are born at the dying point. The process continues in the same way
with independent evolution of each one of the newborn particles. At time $t$
the solution is obtained as a functional of the $n$ existing particles at
time $t$, namely as the product of the initial condition propagated from the
point where each one of the $n$ particles is at time $t$ up to the initial
position. 
\begin{equation}
u(t,x)={\Bbb E}_{x}\left( \varphi _{1}\varphi _{2}\cdots \varphi _{n}\right) 
\label{2.8}
\end{equation}
with 
\begin{eqnarray}
\varphi _{i} &=&\int dy_{1}^{\left( i\right) }dy_{2}^{\left( i\right)
}\cdots dy_{k-1}^{\left( i\right) }dy_{k}^{\left( i\right) }G_{\alpha
,\alpha }^{\beta }\left( \tau _{1},x-y_{1}^{(i)}\right) G_{\alpha ,\alpha
}^{\beta }\left( \tau _{2},y_{1}^{(i)}-y_{2}^{(i)}\right) \cdots   \nonumber
\\
&&\cdots G_{\alpha ,\alpha }^{\beta }\left( \tau
_{k-1},y_{k-2}^{(i)}-y_{k-1}^{(i)}\right) G_{\alpha ,1}^{\beta }\left( \tau
_{k},y_{k-1}^{(i)}-y_{k}^{(i)}\right) u\left( 0^{+},y_{k}^{(i)}\right) 
\label{2.9}
\end{eqnarray}
with $\sum_{j=1}^{k}\tau _{j}=t$, $k-1$ being the number of branchings
leading to particle $i$. Notice that the last propagator in (\ref{2.9}) is
different from the others.

Because of the normalization of the probabilities in the process $B_{\alpha
} $, the probability of each one of the products in (\ref{2.8}) corresponds
to the weight of the corresponding term in the Picard series. Therefore the
expectation value exists whenever the Picard series converges.

The solution (\ref{2.8}) is not yet a purely stochastic solution because it
involves both the expectation value over the process $B_{\alpha }$ and a
multiple integration of the initial condition with the propagation kernels $%
G_{\alpha ,1}^{\beta }$ and $G_{\alpha ,\alpha }^{\beta }$. To obtain a
purely stochastic solution we notice that, for $0<\alpha \leq 1$, the
propagation kernels satisfy the conditions to be the Green's functions of
stochastic processes in ${\Bbb R}$ (see the Appendix).

We denote the processes associated to $G_{\alpha ,1}^{\beta }\left(
t,x\right) $ and $G_{\alpha ,\alpha }^{\beta }\left( t,x\right) $,
respectively by $\Pi _{\alpha ,1}^{\beta }$ and $\Pi _{\alpha ,\alpha
}^{\beta }$. Therefore the process leading to the solution is as described
before with all the particles until the last branching propagating according
to the process $\Pi _{\alpha ,\alpha }^{\beta }$ and the last ones (that
sample the initial condition) propagating by the process $\Pi _{\alpha
,1}^{\beta }$. When finally all the $n$ surviving particles reach time zero,
their coordinates $x+\xi _{i}$ are recorded and the solution is given by 
\begin{equation}
u(t,x)={\Bbb E}_{x}\left( u(0^{+},x+\xi _{1})u(0^{+},x+\xi _{2})\cdots
u(0^{+},x+\xi _{n})\right)   \label{2.11}
\end{equation}
Eq.(\ref{2.11}) is a stochastic solution of (\ref{2.1}) and our main result
is summarized as follows:

{\bf Theorem:} {\it The nonlinear fractional partial differential equation (%
\ref{2.1}), with }$0<\alpha \leq 1${\it , has a stochastic solution given by
(\ref{2.11}), the coordinates }$x+\xi _{i}${\it \ in the arguments of the
initial condition obtained from the exit values of a propagation and
branching process, the branching being ruled by the process }$B_{\alpha }$%
{\it \ and the propagation by }$\Pi _{\alpha ,1}^{\beta }${\it \ for the
particles that reach time }$t${\it \ and by }$\Pi _{\alpha ,\alpha }^{\beta
} ${\it \ for all the remaining ones.}

{\it A sufficient condition for the existence of the solution is} 
\begin{equation}
\left| u(0^{+},x)\right| \leq 1  \label{2.12}
\end{equation}

{\bf Remarks:}

1) The condition $\left| u(0^{+},x)\right| \leq 1$ imposes a finite value
for all contributions to the multiplicative functional. However, the
solution may exist under more general conditions, namely when the decreasing
value of the probability of higher order products in (\ref{2.11})
compensates the growth of the powers of the initial condition.

2) The stochastic solution may also be constructed by a backwards-in-time
stochastic process from time $t$ to time zero. This is obtained by rewriting
Eq.(\ref{2.6}) as 
\begin{eqnarray}
u\left( t,x\right) &=&E_{\alpha ,1}\left( -t^{\alpha }\right) \int_{-\infty
}^{\infty }dy{\cal F}^{-1}\left( \frac{E_{\alpha ,1}\left( -\left( 1+\frac{1%
}{2}\psi _{\beta }^{\theta }\left( k\right) \right) t^{\alpha }\right) }{%
E_{\alpha ,1}\left( -t^{\alpha }\right) }\right) \left( x-y\right) u\left(
0^{+},y\right)  \nonumber \\
&&+\int_{0}^{t}d\tau \tau ^{\alpha -1}E_{\alpha ,\alpha }\left( -\tau
^{\alpha }\right)  \nonumber \\
&&\int_{-\infty }^{\infty }dy{\cal F}^{-1}\left( \frac{E_{\alpha ,\alpha
}\left( -\left( 1+\frac{1}{2}\psi _{\beta }^{\theta }\left( k\right) \right)
\tau ^{\alpha }\right) }{E_{\alpha ,\alpha }\left( -\tau ^{\alpha }\right) }%
\right) \left( x-y\right) u^{2}\left( t-\tau ,y\right)  \label{2.10}
\end{eqnarray}
and noticing that also 
\[
E_{\alpha ,1}\left( -t^{\alpha }\right) +\int_{0}^{t}d\tau \tau ^{\alpha
-1}E_{\alpha ,\alpha }\left( -\tau ^{\alpha }\right) =1 
\]
Then, we obtain the following stochastic construction of the solution:

Starting at time $t$ a particle propagates backwards in time according to
the process $\Pi _{\alpha ,1}^{\beta }$ if it reaches time zero or according
to $\Pi _{\alpha ,\alpha }^{\beta }$ if it branches at time $t-\tau $. The
branching probability is controlled by the process $B_{\alpha }$ (that is,
the branching probability density is $\tau ^{\alpha -1}E_{\alpha ,\alpha
}\left( -\tau ^{\alpha }\right) $). When it branches, two new particles are
born which propagate independently and the process is repeated until all
surviving particles reach time zero.\bigskip

{\LARGE Appendix. The Green's functions and the characterization of the
processes}

{\bf The processes }$\Pi _{\alpha ,1}^{\beta }$ {\bf and} $\Pi _{\alpha
,\alpha }^{\beta }$%
\begin{equation}
{\cal F}\left\{ G_{\alpha ,1}^{\beta }\left( t,x\right) \right\} \left(
t,k\right) =\frac{E_{\alpha ,1}\left( -\left( 1+\frac{1}{2}\psi _{\beta
}^{\theta }\left( k\right) \right) t^{\alpha }\right) }{E_{\alpha ,1}\left(
-t^{\alpha }\right) }  \label{A.1}
\end{equation}
\begin{equation}
{\cal F}\left\{ G_{\alpha ,\alpha }^{\beta }\left( t,x\right) \right\}
\left( t,k\right) =\frac{E_{\alpha ,\alpha }\left( -\left( 1+\frac{1}{2}\psi
_{\beta }^{\theta }\left( k\right) \right) t^{\alpha }\right) }{E_{\alpha
,\alpha }\left( -t^{\alpha }\right) }  \label{A.2}
\end{equation}

For a propagation kernel $G\left( t,x\right) $ to be the Green's function of
a stochastic process, the following conditions should be satisfied:

(i) $G\left( 0,x-y\right) =\delta \left( x-y\right) $ or ${\cal F}\left\{
G\right\} \left( 0,k\right) =1$ $\forall k$

(ii) $\int dxG\left( t,x\right) =1$ $\forall t$ or ${\cal F}\left\{
G\right\} \left( t,0\right) =1$

(iii) $G\left( t,x\right) $ should be real and $\geq 0$

For the processes $\Pi _{\alpha ,1}^{\beta }$ and $\Pi _{\alpha ,\alpha
}^{\beta }$

(i) ${\cal F}\left\{ G_{\alpha ,1}^{\beta }\right\} \left( 0,k\right) =\frac{%
E_{\alpha ,1}\left( 0\right) }{E_{\alpha ,1}\left( 0\right) }=1$ and ${\cal F%
}\left\{ G_{\alpha ,\alpha }^{\beta }\right\} \left( 0,k\right) =\frac{%
E_{\alpha ,\alpha }\left( 0\right) }{E_{\alpha ,\alpha }\left( 0\right) }=1$

(ii) ${\cal F}\left\{ G_{\alpha ,1}^{\beta }\right\} \left( t,0\right) =%
\frac{E_{\alpha ,1}\left( -t^{\alpha }\right) }{E_{\alpha ,1}\left(
-t^{\alpha }\right) }=1$ and ${\cal F}\left\{ G_{\alpha ,\alpha }^{\beta
}\right\} \left( t,0\right) =\frac{E_{\alpha ,\alpha }\left( -t^{\alpha
}\right) }{E_{\alpha ,\alpha }\left( -t^{\alpha }\right) }=1$

(iii) If ${\cal F}\left\{ G\right\} \left( t,-k\right) =\left( {\cal F}%
\left\{ G\right\} \left( t,k\right) \right) ^{*}$ then $G\left( t,x\right) $
is real.

Because $\psi _{\beta }^{\theta }\left( -k\right) =\left( \psi _{\beta
}^{\theta }\left( k\right) \right) ^{*}$ it follows 
\[
E_{\alpha ,1}\left( -\left( 1+\frac{1}{2}\psi _{\beta }^{\theta }\left(
-k\right) \right) t^{\alpha }\right) =\left( E_{\alpha ,1}\left( -\left( 1+%
\frac{1}{2}\psi _{\beta }^{\theta }\left( k\right) \right) t^{\alpha
}\right) \right) ^{*} 
\]
and 
\[
E_{\alpha ,\alpha }\left( -\left( 1+\frac{1}{2}\psi _{\beta }^{\theta
}\left( -k\right) \right) t^{\alpha }\right) =\left( E_{\alpha ,1}\left(
-\left( 1+\frac{1}{2}\psi _{\beta }^{\theta }\left( k\right) \right)
t^{\alpha }\right) \right) ^{*} 
\]
implying that both $G_{\alpha ,1}^{\beta }\left( t,x\right) $ and $G_{\alpha
,\alpha }^{\beta }\left( t,x\right) $ are real.

Finally, for the positivity, one notices that for $0<\alpha \leq 1$ and $%
\rho \geq \alpha $, $E_{\alpha ,\rho }\left( -x\right) $ is a completely
monotone function\cite{Schneider}. Therefore 
\[
E_{\alpha ,\rho }\left( -x\right) =\int_{0}^{\infty }e^{-rx}dF\left(
r\right) 
\]
with $F$ nondecreasing and bounded.

For $G_{\alpha ,\rho }^{\beta }\left( t,x\right) $ ($\rho =1$ and $\rho
=\alpha $) one has 
\begin{eqnarray*}
G_{\alpha ,\rho }^{\beta }\left( t,x\right) &=&\frac{1}{2\pi E_{\alpha ,\rho
}\left( -t^{\alpha }\right) }\int_{0}^{\infty }dF\left( r\right)
\int_{-\infty }^{\infty }dke^{-ikx}e^{-rt^{\alpha }\left( 1+\frac{1}{2}\psi
_{\beta }^{\theta }\left( -k\right) \right) } \\
&=&\frac{1}{2\pi E_{\alpha ,\rho }\left( -t^{\alpha }\right) }%
\int_{0}^{\infty }dF\left( r\right) e^{-rt^{\alpha }}\int_{-\infty }^{\infty
}dke^{-ikx}e^{-\frac{rt^{\alpha }}{2}\psi _{\beta }^{\theta }\left(
-k\right) }
\end{eqnarray*}
We recognize the last integral (in $k$) as the Green's function of a Levy
process. Therefore one has an integral in $r$ of positive quantities
implying that $G_{\alpha ,1}^{\beta }\left( t,x\right) $ and $G_{\alpha
,\alpha }^{\beta }\left( t,x\right) $ are positive.

{\bf The process }$B_{\alpha }$

The decaying probability in time $d\tau $ of this process is 
\[
\tau ^{\alpha -1}E_{\alpha ,\alpha }\left( -\tau ^{\alpha }\right) 
\]
From 
\[
\int_{0}^{t}\tau ^{\alpha -1}E_{\alpha ,\alpha }\left( -\tau ^{\alpha
}\right) d\tau =1-E_{\alpha ,1}\left( -t^{\alpha }\right) 
\]
it follows that $E_{\alpha ,1}\left( -t^{\alpha }\right) $ is the survival
probability up to time $t$. The process $B_{\alpha }$ is a fractional
generalization of the exponential process.

\end{document}